\documentclass{amsart}

\usepackage{amsmath,graphicx, hyperref}

\newtheorem{theorem}{Theorem}[section]
\newtheorem{proposition}[theorem]{Proposition}
\newtheorem{lemma}[theorem]{Lemma}
\newtheorem{corollary}[theorem]{Corollary}
\theoremstyle{definition}
\newtheorem{definition}[theorem]{Definition}
\newtheorem{remark}[theorem]{Remark}
\newtheorem{convention}[theorem]{Convention}

\newcommand{\GammaH}{\Gamma (G, X\cup \mc H) }
\newcommand{\mc}[1]{\mathcal{#1}}
\newcommand{\N}{\ensuremath{\mathbb{N}}}
\DeclareMathOperator{\area}{\mathrm{Area}}
\DeclareMathOperator{\dist}{\mathsf{dist}}
\DeclareMathOperator{\Label}{\mathrm{Label}}

\begin{document}

\title[Quasiconvexity and Relative Hyperbolicity]{On Quasiconvexity and Relatively Hyperbolic Structures on Groups}
\author[E.Mart\'inez-Pedroza]{Eduardo Mart\'inez Pedroza}
      \address{McMaster University\\
               Hamilton, Ontario, Canada L8P 3E9}
      \email{eduardo.martinez.pedroza@gmail.com}

\begin{abstract}
Let $G$ be a group which is hyperbolic relative to a collection of subgroups $\mc{H}_1$, and it is also hyperbolic relative to a collection of subgroups $\mc{H}_2$.
Suppose that $\mc{H}_1 \subset \mc{H}_2$. We characterize when a relative quasiconvex subgroup of $(G, \mc H_2)$ is still relatively quasiconvex in $(G, \mc H_1)$.  We also show that relative quasiconvexity is preserved when passing from $(G, \mc H_1)$ to $(G, \mc H_2)$.  Applications are discussed.
\end{abstract}

\maketitle

\section{Introduction}\label{sec.introduction}

If $G$ is a countable group and $\mc H$ is a finite collection of subgroups of $G$, the notion of (strong) relative hyperbolicity for the pair $(G, \mc H)$  has been defined by different authors,  all these definitions being equivalent~\cite{HK08, Os06}.  Equivalent notions of relatively quasiconvex subgroups for relatively hyperbolic groups have also been introduced by different authors and they have been proved to be equivalent, we refer the reader to~\cite{HK08, MaWi10b}.

The collection of quasiconvex subgroups of a (strong)  relatively hyperbolic pair $(G, \mc H)$ is an interesting class of subgroups of $G$. For example, this class is closed under finite intersections~\cite{HK08, MP09}, contains all virtually cyclic subgroups~\cite{Os06}, and every element naturally admits a  relatively hyperbolic group structure~\cite{HK08, MaWi10b}. A countable group $G$ can have different relatively hyperbolic structures giving rise to different collections of subgroups of $G$ with these properties.  In this paper, we start an investigation of how the collection of quasiconvex subgroups of $(G, \mc H)$ varies with respect to $\mc H$. 
Our main result is the following:

\begin{theorem}\label{thm:main}
Let $G$ be a countable  group, and let $(G, \mc H_1)$ and $(G, \mc H_2)$ be relatively hyperbolic structures with $\mc H_1 \subset \mc H_2$ and $\mc H_2$ finite. 
\begin{enumerate}
\item \label{main-1} If $Q$ is a quasiconvex subgroup of $(G, \mc H_1)$, then $Q$ is quasiconvex in $(G, \mc H_2)$.
\item \label{main-2} If $Q$ is a quasiconvex subgroup of $(G, \mc H_2)$, and for each $H \in \mc H_2\setminus \mc H_1$ and each $g \in G$ the subgroup $Q\cap gHg^{-1}$ is quasiconvex in $(G, \mc H_1)$,  then $Q$ is quasiconvex in $(G, \mc H_1)$.
\end{enumerate}
\end{theorem}

\begin{remark}
During the review process of this work, Wen-Yuan Yang extended the results of this paper in~\cite{Ya11}. In particular, Theorem~\ref{thm:main} holds for a more general extension of finite peripheral structures: for each subgroup $H\in \mathcal{H}_1$ there is a subgroup $H'\in \mathcal{H}_2$ such that $H<H'$. 
\end{remark}

The proofs of the theorem and the corollaries stated below are discussed in Section~\ref{sec:proofs}. In the corollaries, all relatively hyperbolic pairs $(G, \mc H)$ consist of a countable group $G$ and a finite collection of subgroups $\mc H$.

\subsection*{Elementary subgroups in peripheral structures.}
Recall that a group is called \emph{elementary} if it contains a cyclic group of finite index. 
\begin{corollary}\label{cor:elementary}
Let $G$ be a  countable group, and $(G, \mc H_1)$ and $(G, \mc H_2)$ be relatively hyperbolic structures with $\mc H_1 \subset \mc H_2$.
Then $(G, \mc H_1)$ and $(G, \mc H_2)$ have the same class of quasiconvex subgroups if and only if every $H \in \mc H_2 \setminus \mc H_1$ is elementary.
\end{corollary}

\subsection*{Coherence and Local Quasiconvexity}
\begin{definition}
A group is \emph{coherent} if every finitely generated subgroup is finitely presented.  A relatively hyperbolic pair $(G, \mc H)$ is called \emph{locally quasiconvex} if every finitely generated subgroup of $G$ is a quasiconvex subgroup of $(G, \mc H)$.  
\end{definition}
The following result follows from the fact that relative quasiconvex subgroups naturally admit relative hyperbolic structures~\cite{HK08, MaWi10b} together with results of~\cite{Os06}, see Section~\ref{sec:proofs} for a proof.
\begin{theorem}\label{thm:coh}
Suppose that $(G, \mc H)$ is relatively hyperbolic and locally quasiconvex. If each $H \in \mc H$ is a coherent group, then $G$ is coherent. 
\end{theorem}
Using Theorem~\ref{thm:main},  the following  result analogous to Theorem~\ref{thm:coh} is obtained.
\begin{corollary}\label{cor:coh}
Suppose that $(G, \mc H)$ is relatively hyperbolic and locally quasiconvex. 
If each $H \in \mc H$ is a locally quasiconvex hyperbolic group, then $G$ is a locally quasiconvex hyperbolic group.
\end{corollary}
\begin{remark} In~\cite{MaWi10}, Wise and the author use Theorem~\ref{thm:coh} and Corollary~\ref{cor:coh} to obtain results on coherence and local quasiconvexity of high powered one relator products.
\end{remark}

\subsection*{Malnormal Cores of Quasiconvex Subgroups}
\begin{definition}[Essentially Distinct Conjugates]
Let $Q$ be a subgroup of a group $G$. Let  $\{g_i | 1\leq i \leq n\}$ be a collection of elements of $G$. 
Conjugates $\{g_iQg_i^{-1} | 1\leq i\leq n\}$ are called \emph{essentially distinct} if $g_iQ \neq g_j Q$ for $i\neq j$.
\end{definition}

We recover the following result from~\cite[Proposition 3.12]{AGM08} without the assumption that the group $G$ is torsion free. 
\begin{corollary} \label{cor:3}
Let $Q$ be a quasiconvex subgroup of a hyperbolic group $G$. Then there is a finite collection subgroups $\mc{H}$ of $G$ with the following properties.
\begin{enumerate}
\item\label{cor:3-1} $(G , \mc H)$ is relatively hyperbolic. 
\item\label{cor:3-2} $Q$ is quasiconvex in $(G , \mc H)$.  
\item\label{cor:3-3}  For each $H \in \mc H$, the subgroup $Q\cap H$ has finite index in $H$.
\item\label{cor:3-4}  If $\{g_iQg_i^{-1} | 1\leq i\leq n\}$ is a maximal collection of essentially distinct conjugates of $Q$ with infinite intersection.
Then $\bigcap g_iQg_i^{-1}$ is a (finite index) subgroup of some $H \in \mc H$ up to conjugacy in $G$.
\end{enumerate}
\end{corollary}

\begin{remark}
The construction of $\mc H$ in Corollary~\ref{cor:3}  is based on the fact quasiconvex subgroups have finite height, and have finite index in their commensurators~\cite{GMRS, KS96}. Statements~\eqref{cor:3-1},~\eqref{cor:3-3} and~\eqref{cor:3-4} of Corollary~\ref{cor:3} follow directly from the construction of $\mc H$. A non-trivial argument for statement~\eqref{cor:3-2} is given in~\cite{AGM08} assuming that $G$ is torsion free.  This is an immediate consequence of Theorem~\ref{thm:main} without the assumption on torsion freeness. 
\end{remark}

\subsection*{Combination of Hyperbolically Embedded Subgroups}

\begin{definition}\cite{Os06-3}\label{def.hyp_embedded}
Let $(G, \mc H)$ be relatively hyperbolic.  A subgroup $P \leq G$ is said to be {\it hyperbolically embedded into $(G,  \mc H)$}, if the pair $(G, \mc{H} \cup \{P\})$ is relatively hyperbolic. 
\end{definition}

\begin{corollary} \label{cor:1}
Let $(G, \mc H)$ be relatively hyperbolic, and suppose that $X$ is a finite generating set of $G$.
For every quasiconvex subgroup $Q$ of $(G, \mc H)$, and every hyperbolically embedded subgroup $P$ of $(G, \mc H)$, 
there is constant $C > 0$ with the following property. 
If $R$ is a subgroup of $P$ such that
\begin{enumerate}
\item $Q \cap P \leq R$, and
\item $|g|_X \geq C$ for any $g \in R \setminus Q$,
\end{enumerate}
then the natural homomorphism $Q \ast_{Q\cap R} R \longrightarrow G$
is injective. Moreover, if the subgroup $R$ is quasiconvex relative to $\mathcal{H}$, then the subgroup $\langle Q \cup R\rangle$ is quasiconvex relative to $\mathcal{H}$.
\end{corollary}
\begin{remark}
Corollary~\ref{cor:1} follows from Theorem~\ref{thm:main} and~\cite[Theorem 1.1]{MP09}. In the case the $G$ is a hyperbolic group, Corollary~\ref{cor:1} is a result of Gitik~\cite{Gi99}.
\end{remark}

\begin{corollary}\label{cor:2}
Let $G$ be a non-elementary and properly relatively hyperbolic group with respect to a collection of subgroups $\mathcal{H}$.
Suppose that $G$ has no non-trivial finite normal subgroups. Then for any finite subset $F$ of non-trivial elements of $G$ there exists an element $g \in G$ with the following properties. For every $f \in F$,
\begin{enumerate}
\item $\langle f, g \rangle$ is isomorphic to the free product $\langle f \rangle \ast \langle g \rangle$, and
\item $\langle f, g \rangle$ is a quasiconvex subgroup relative to $\mathcal{H}$.
\end{enumerate}
\end{corollary}
\begin{remark}
Corollary~\ref{cor:2}   is a stronger version of a result by G. Arzhntseva and A. Minasyan~\cite[Theorem 1]{AM07} used to prove that relatively hyperbolic groups with no non-trivial finite normal subgroups are $C^*$-simple. The original result does not include the statement on quasiconvexity.
\end{remark}

\subsection*{Outline of the Paper.} 

 Section~\ref{sec.hyperbolicity} recalls the definition of relative hyperbolic groups and some geometric results  of  Osin in~\cite{Os06}. Then we extend a result on the geometry of finitely generated relatively hyperbolic groups~\cite[Theorem 3.23]{Os06} for countable (not necessarily finitely generated) relatively hyperbolic groups, Theorems~\ref{thm:BCP} and~\ref{thm:BCP-countable}.

Section~\ref{sec.quasiconvexity} recalls a definition of relative quasiconvexity by Hruska in~\cite{HK08} for countable relatively hyperbolic groups. Then we provide short and direct proofs of some known results for relatively quasiconvex subgroups using Theorem~\ref{thm:BCP-countable}.  

Section~\ref{sec.quasigeodesics} recalls a characterization of hyperbolically embedded subgroups by Osin in~\cite{Os06-3}. Then we use this characterization
together with some results of Section~\ref{sec.hyperbolicity} to obtain  a result on the geometry of  hyperbolically embedded subgroups stated as Proposition~\ref{prof:hyp_emb}. 

Section~\ref{sec:proofs} contains the proofs of Theorem~\ref{thm:main} and the corollaries.

\subsection*{Acknowledgements:}  The author thanks the referee for useful comments, suggestions and critical corrections.  
The author thanks Noel Brady, Chris Hruska, Jason Manning, Ashot Minasyan for comments on this work.
This project was supported by a Britton Postdoctoral Fellowship at McMaster University.

\section{Relative Hyperbolicity} \label{sec.hyperbolicity}

This section consists of three parts. The first and second part recall the definition of relative hyperbolicity and some results  from~\cite{Os06}. In the third part, we revisit some results of~\cite[Chapter 3]{Os06} that require the assumption of finite generation.  We state these more general versions as Theorems~\ref{thm:BCP} and~\ref{thm:BCP-countable} and give accounts of the modifications that are required in their proofs.

\subsection{Definition of Relative hyperbolicity by Osin}\label{sec:def} The following definition is taken from~\cite[Chapter 2]{Os06}.

Let $G$ be a group, $\{ H_\lambda \}_{\lambda \in \Lambda}  $ be a collection of subgroups of $G$, and $X$ a subset of $G$.
The collection $X$ is a \emph{relative generating set of $G$ with respect to $\{ H_\lambda \}_{\Lambda}  $} if  $G$ is generated by $X\cup \bigcup_{\lambda \in \Lambda} H_\lambda$ and $X$ is closed under inverses.  In this case we say that \emph{$X$ is a relative generating set of $(G, \{ H_\lambda \}_{ \Lambda})$}.

Suppose that $X$ is relative generating set with respect to $\{ H_\lambda \}_{\Lambda}  $. 
In this case, the group $G$ is a quotient of the free product $F$
\begin{equation}\label{eq:F} F = ( \ast_{\lambda \in \Lambda} H_\lambda ) \ast F(X) \end{equation}
where $F(X)$ denotes the free group with basis $X$.  Let $N$ be the kernel of the natural quotient map $F\rightarrow G$,
and suppose that $N$ is the normal closure of a subset $\mc R$ of $F$. Then we say that $G$ has a \emph{relative presentation}
\begin{equation}\label{eq:rel-pr} \langle  X,\  H_\lambda, \lambda \in \Lambda \  |  \  \mc R \rangle \end{equation}
If $X$ and $\mc R$ are finite sets, the relative presentation~\eqref{eq:rel-pr} is called \emph{finite}, and $G$ is called \emph{finitely presented relative to $\{ H_\lambda \}_{\Lambda}  $}.

Let
\[ \mc H = \bigsqcup_{\lambda \in \Lambda} (  H_\lambda \setminus \{1\} ).\]
Given a word $W$ in the free monoid $(X\cup \mc H)^*$ representing the identity in $G$, there exists an expression
\begin{equation}\label{eq:kernel} W =_F  \prod_{i=1}^k f_i^{-1}R_if_i \end{equation}
with equality in the group $F$ defined in~\eqref{eq:F}, and where $R_i \in \mc R$, and $f_i \in F$ for each $1\leq i \leq k$.
The smallest possible value of $k$ for an expression representing $W$ as in~\eqref{eq:kernel} is denoted by $\area^{rel}(W)$. 
A function $f:\N \rightarrow \N$ is a \emph{relative isoperimetric function} of~\eqref{eq:rel-pr} if $\area^{rel}(W) \leq f(||W||)$  for every word $W$.

\begin{definition}[Relative Hyperbolicity]\label{def:rel-hyp}
A group \emph{$G$ is hyperbolic relative to a collection of subgroups $\{ H_\lambda \}_{\Lambda}  $} if there is a finite relative
presentation of $G$ with respect to $\{ H_\lambda \}_{\Lambda}  $ admitting a linear relative isoperimetric function.
In this case, we say that the pair $(G, \{ H_\lambda \}_\Lambda)$ is a \emph{relatively hyperbolic structure}.
\end{definition}

\begin{theorem}\cite[Theorem  2.34]{Os06}
Definition~\ref{def:rel-hyp} is independent of the choice of finite relative presentation. 
\end{theorem}

\begin{theorem}\cite[Theorem  1.4]{Os06} \label{thm:malnormality}
Let $G$	 be hyperbolic relative to a collection of subgroups $\{ H_\lambda  \}_{\lambda \in \Lambda}$. Then
\begin{enumerate}
\item For any $g_1,g_2 \in G$, the intersection $H_\lambda^{g_1} \cap H_\mu^{g_2}$ is finite whenever $\lambda \neq \mu$.
\item The intersection $H_\lambda^{g} \cap H_\lambda$ is finite for any $g \not \in H_\lambda$.
\end{enumerate}
\end{theorem}

\subsection{Geometry of Relatively Hyperbolic Groups, the Main Tools.}

Let $G$ be hyperbolic relative to a collection of subgroups $\{ H_\lambda \}_{\Lambda}  $,  let $X$  be a finite relative generating set of $G$ with respect to $\{ H_\lambda \}_{\Lambda}  $, and let $\mc H = \bigsqcup_{\lambda \in \Lambda} H_\lambda \setminus \{1\}$.

In this part, we recall some terminology and results of~\cite[Chapter 2 ]{Os06} on the geometry of relatively hyperbolic groups.  

\begin{definition}[Relative Cayley Graph]
Recall that the \emph{Cayley graph $\Gamma (G, X)$ of $G$ with respect to the subset $X \subset G$} is the graph with vertex set $V(\Gamma) =G$,
edge set $E(\Gamma)= G \times X$, and an edge $e=(g, X)$ goes from $g$ to $gs$.  The graph $\Gamma (G, X)$ is connected  if and only if $X$
is a generating set of $G$.  Let $\Gamma(G, X\cup \mc H)$ denote  the Cayley graph of $G$ with respect to the generating set $X\cup \mc H$. 
\end{definition}

\begin{theorem} \cite[Corollary 2.54]{Os06} \label{thm:hyperbolicity}
Let $G$ be a group and let $\{ H_\lambda \}_{\Lambda}  $ be a collection of subgroups.
Suppose that there is a finite relative presentation of $G$ with respect to $\{ H_\lambda \}_{\Lambda}  $ admitting 
a relative isoperimetric function.  Then the following statements are equivalent.
\begin{enumerate}
\item $G$ is hyperbolic relative to $\{ H_\lambda \}_{\Lambda}  $.
\item  The Cayley graph $\Gamma (G, X\cup \mc H)$ is a $\delta$-hyperbolic metric space for some $\delta \geq 0$.
\end{enumerate}
\end{theorem}

\begin{definition}[Basic terminology]
Let $p=e_1\cdots e_k$ be a combinatorial path in the Cayley graph $\Gamma(G, X\cup \mc H)$.  We refer to its initial vertex as $p_-$, and its terminal vertex as $p_+$.  The path $p$ determines a word \[ \Label(p) = \Label (e_1) \Label (e_2) \cdots \Label (e_k)\] in the alphabet $X\cup \mc H$ that represents an element $g \in G$ so that $p_+ =  p_- g$. The \emph{length} of $p$ is denoted by  $\ell (p)$. 
\end{definition}

\begin{definition}[$\mc H$--components, Phase vertices]
Sub-paths of $q$ with at least one edge are called \emph{non-trivial}.  An \emph{$H_\lambda$--component} of $q$ is a maximal non-trivial sub-path $X$ of $q$ with $\Label(X)$ a word in the (sub)-alphabet $H_\lambda$. When it is not necessary to specify the subgroup $H_\lambda$, we refer to the $H_\lambda$--component as
an $\mc H$--component.  Two $H_\lambda$--components $s_1$, $s_2$ are \emph{connected} if the vertices of $s_1$ and $s_2$ belong to the same left coset of $H_\lambda$.    An $H_\lambda$--component $X$ of $q$ is \emph{isolated} if it is not connected to a different $H_\lambda$--component of $q$.  A vertex $v$ of $q$ is called \emph{ phase} if it is not an interior vertex of a $\mc H$--component $X$ of $q$.  
\end{definition}

\begin{convention}
For $X \subset G$ and $g\in G$, $|g|_X$ denotes the word length with respect to the set $X$ if $g \in \langle X \rangle$, otherwise $|g|_X$ is defined as infinity.  In particular, each $X \subset G$ defines a metric $\dist_X$ with values in $\mathbb{R}\cup \{ \infty \}$ on the vertices of $\Gamma (G, X \cup \mc H)$.
\end{convention}

Proposition~\ref{thm:n-gon-2} is a powerful technical result due to Osin~\cite{Os06, Os06-3} on which the results of this paper rely on. 
\begin{proposition} \cite[Lemma 2.11]{Os06-3} \label{thm:n-gon-2}
There exists a constant $M>0$ and subsets $\Omega_\lambda \subset H_\lambda$ for each $\lambda \in \Lambda$ such that the following hold.
\begin{enumerate}
\item  The union $\Omega = \bigcup_{\lambda \in \Lambda} \Omega_\lambda$ is finite.
\item  Let $q$ be a cycle in $\Gamma (G, X\cup \mc H)$, let $p_1, \cdots , p_k$ be a collection of isolated $H_\lambda$-components of $q$ for some $\lambda \in \Lambda$,  and let $g_1, \dots , g_k$ be the elements of $G$ represented by the labels of $p_1, \cdots , p_k$.  Then for any $i=1, \dots , k$, the element $g_i$ belongs
to the subgroup $\langle \Omega_\lambda \rangle$ and the lengths of $g_i$ with respect to $\Omega_\lambda$ satisfy the inequality
\[ \sum_{i=1}^k |g_i|_{\Omega_\lambda} \leq M \cdot \ell (q) .\]
\end{enumerate}
\end{proposition}

\subsection{The Bounded Coset Penetration Property}
In this subsection, we revisit~\cite[Theorem. 3.23]{Os06} which is a result for finitely generated relatively hyperbolic groups known as the Bounded Coset Penetration property. This property was originally introduced by Farb in~\cite{Fa98} as part of an equivalent definition of relative hyperbolicity in the finite generating case. 

The main results of this section are Theorems~\ref{thm:BCP} and~\ref{thm:BCP-countable} which are versions of~\cite[Theorem. 3.23]{Os06} for general
relatively hyperbolic groups and for countable relatively hyperbolic groups respectively. The  proof of~\cite[Theorem. 3.23]{Os06} relies on~\cite[Lemma 3.1]{Os06} which is replaced by the more general result~\cite[Lemma 2.27]{Os06} stated in this paper as Proposition~\ref{thm:n-gon-2}. This is the most important modification in the proof besides some extra estimations and minor observations. Below there is an account of these modifications of Osin's arguments in~\cite{Os06}. These modifications have been also pointed out in the work of Yang~\cite{Ya11}.

Let $G$ be hyperbolic relative to a collection of subgroups $\{ H_\lambda \}_{\Lambda} $,  let  $X$  be a finite relative generating set of $G$ with respect to $\{ H_\lambda \}_{\Lambda}  $, let $\mc H = \bigcup_{\lambda \in \Lambda} H_\lambda \setminus \{1\}$, and let $\Omega$ be the finite subset of $G$ provided by
Proposition~\ref{thm:n-gon-2}.

\begin{definition}[$K$-similar paths] \label{def:similar}
For a $k>0$, two paths $p$ and $q$ in $\Gamma (G, X\cup \mc H)$ are $k$-similar if $ \max\{ \dist_{X} (p_-, q_-),  \dist_{X}(p_+, q_+) \}  \leq k .$
\end{definition}

\begin{definition}[Backtracking]
A path $p$ in $\Gamma (G, X\cup \mc H)$ is a \emph{path without backtracking} if every $\mathcal{H}$--component of $p$ is isolated.
\end{definition}

\begin{definition}[Quasi-geodesic]
A path $q$ in a metric space $(Y, \dist)$ is a \emph{$(\mu, c)$-quasi-geodesic path} for some $\mu \geq 1$, $c\geq 0$, if for every subpath $p$ of $q$ 
the inequality $ \ell (p) \leq \mu \dist (p_-, p_+) +c$ holds, where $\ell (p)$ denotes the length of the subpath $p$. A \emph{geodesic path} is a $(1,0)$-quasi-geodesic.
\end{definition}
\begin{remark}
A geodesic path $q$ in $\Gamma(G, X\cup \mc H)$ has not backtracking, all $\mathcal{H}$--components of $q$ consist of a single edge, and all vertices of $q$ are phase. 
\end{remark}

\begin{theorem}[Bounded Coset Penetration Property]\label{thm:BCP}
For each $\mu \geq 1, c \geq 0, k\geq 0$, there exists $\epsilon=\epsilon (\mu, c, k) > 0$ with the following property. 
For any two $k$-similar $(\mu, c)$-quasi-geodesics without backtracking $p$ and $q$ in $\Gamma (G, X\cup \mc H)$, the following conditions hold.
\begin{enumerate}
\item The sets of phase vertices of the paths $p$ and $q$ are contained in the closed $\epsilon$-neighborhoods with respect to the metric $\dist_{X\cup \Omega}$ of each other.
\item If $s$ is an $H_\lambda$--component of $p$ such that $\dist_{X\cup \Omega} (s_-, s_+) \geq \epsilon $, then there exists an $H_\lambda$--component   $t$ of $q$ which is connected to $s$.
\item If $s$ and $t$ are connected $H_\lambda$--components of $p$ and $q$ respectively, then 
\[\max\{ \dist_{X\cup \Omega }(s_-, t_-), \dist_{X\cup \Omega }(s_+, t_+) \} \leq \epsilon.\]
\end{enumerate}
\end{theorem}

\begin{remark}
The only difference between the the statements of~\cite[Theorem. 3.23]{Os06} and Theorem~\ref{thm:BCP} consists of replacing in the conclusion of~\cite[Theorem. 3.23]{Os06} the metric $\dist_X$ with $\dist_{X\cup \Omega}$.  This enlargement of the finite generating set is necessary as the following example provided by the referee of this paper illustrates.  If $G$ is an amalgamated product $H_1\ast_C H_2$ over a finite subgroup $C$, then $G$ is hyperbolic relative to $\{H_1, H_2\}$ and we can take the empty set as the relative generating set $X$. In this case,  if $p$ and $q$ are $k$-similar geodesics, the sets of phase vertices of $p$ and $q$ are contained in the closed $\epsilon$-neighborhoods with respect to the metric $\dist_{X}$ of each other if only if the $p$ and $q$ are the same path. However in the case that $C$ is non-trivial, there are different geodesics in $\Gamma (G, X\cup \mc H)$ between the same pair of vertices.  This was also pointed out by Wen-Yuan Yang 
in~\cite[Example 2.12]{Ya11}.
\end{remark}

Theorem~\ref{thm:BCP} follows by taking together Propositions~\ref{prop:BCP-1},~\ref{prop:BCP-2}, and~\ref{prop:BCP-3} below. 
These Propositions are more general versions of~\cite[Propositions 3.15, Lemma 3.21, Lemma 3.22]{Os06} respectively. 
 
\begin{lemma} \cite[Lemma 3.8]{Os06}. \label{lem:similar}
For any $\delta\geq 0$, $\mu \geq 1$, $c\geq 0$, and $k\geq 0$, there exists a constant $K=K(\delta, \mu, c, k)$ such that 
the following holds. Suppose that $(Y,\dist )$ is a $\delta$-hyperbolic space and $p,q$ are $(\mu, c)$-quasigeodesics paths in $Y$
such that \[\max\{ \dist (p_-, q_-), \dist (p_+, q_+) \} \leq k.\] Then $p$ and $q$ belong to the closed $K$-neighborhoods of each other.
\end{lemma}

\begin{figure}
\begin{center} 
\includegraphics[width=0.9\linewidth]{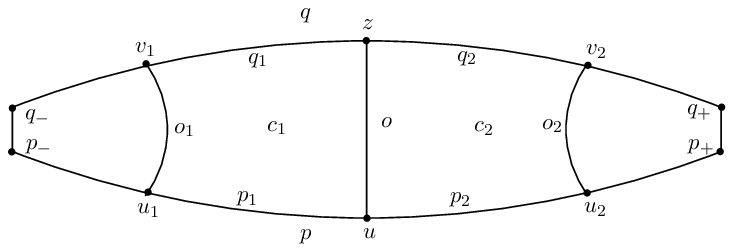}%
\end{center} 
\caption{Paths in the proof of Prop.~\ref{prop:BCP-1}, ~\cite[Fig. 3.1]{Os06}. } 
 \label{fig:Figura-1} 
\end{figure}

\begin{proposition}[Bounded Coset Penetration Property, Part I]\label{prop:BCP-1}
For each $\mu \geq 1, c \geq 0, k\geq 0$, there exists $\epsilon=\epsilon (\mu, c, k) > 0$ with the following property. 
Let  $p$ and $q$ by two $k$-similar $(\mu, c)$--quasi-geodesic paths in $\Gamma(G, X\cup \mc H)$ such that $p$ is a path without backtracking.
Then for any phase vertex $u$ of $p$ there exists a phase vertex $v$ of $q$ such that $ \dist_{X\cup \Omega}(u,v) \leq \epsilon .$
\end{proposition}
Minor modifications in the proof by Osin of~\cite[Proposition 3.15]{Os06} yield a proof of Proposition~\ref{prop:BCP-1}. Specifically, the proof of~\cite[Proposition 3.15]{Os06} is divided into two parts. First the \emph{main geometric construction} where the only extra observation required is that for $a,b\in G$ there are paths in the possibly disconnected Cayley graph $\Gamma (G, X)$ between $a,b$ whenever $\dist_X (a,b) <\infty$. In the second part of the proof, the \emph{final estimation}, one requires to apply Proposition~\ref{thm:n-gon-2} instead of~\cite[Lemma 3.1]{Os06}, and check some variations of the computations to define $\epsilon_1 (\mu, c, k)$. 
\begin{proof}[Proof of Proposition~\ref{prop:BCP-1}] 
First observe that by replacing each $\mc H$-component of $p$ and $q$ by a single edge, one obtains paths which are still a $(\mu, c)$-quasigeodesic and have no backtracking (see~\cite[Lemma 3.12]{Os06}). Hence we can assume that all vertices of $p$ and $q$ are phase. Let $u$ be a vertex of $p$. 

Since $\Gamma(G, X\cup \mc H)$ is a $\delta$-hyperbolic space, let $K_0=K(\delta, \mu, c, k)$  and $K=K(\delta, \mu, c, K_0)+1/2$ be the constants
provided by Lemma~\ref{lem:similar}. Without loss of generality, assume that $K$ and $K_0$ are integers and  that $ K\geq K_0\geq k$. 

{\bf Initial geometric construction:} We reproduce illustration~\cite[Fig 3.1]{Os06} as Figure~\ref{fig:Figura-1}.  Seven paths denoted by $o$ and $p_i, q_i, o_i$ for $i=1,2$, and a vertex $z$ on $q$  are chosen:

\emph{Paths $p_i$.} Let $u_1$ be a vertex on the segment $[p_-, u]$ of $p$ such that either $u_1=p_-$ and $\dist_{X\cup \mc H}(p_-, u) \leq 2K$, or $\dist_{X\cup \mc H}(u_1, u) = 2K$. Analogously, choose $u_2$ on the segment $[u, p_+]$ of $p$ such that either $u_2=p_+$ and $\dist_{X\cup \mc H}(p_-, u) \leq 2K$, or $\dist_{X\cup \mc H}(u, u_2) = 2K$. Let $p_1$ and $p_2$ be the segments $[u_1, u]$ and $[u, u_2]$ of $p$ respectively. Since $p$ is a $(\mu, c)$-quasi-geodesic, 
\begin{equation}\label{eq:bcp-p} \max\{ \ell (p_1) , \ell (p_2) \} \leq 2 \mu K + c. \end{equation}

\emph{Paths $o_i$.}  Choose vertices $v_1,v_2$ of $q$, and two paths $o_1, o_2$  in $\Gamma (G, X\cup \mc H)$ such that $(o_i)_-=u_i$, $(o_i)_+=v_i$, $i=1,2$, and  \begin{equation}\label{eq:bcp-oi}  \ell (o_i) \leq K_0 \leq K,\end{equation}
The choices are made so that: if $u_1=p_-$, then $v_1=q_-$ and  $o_1$ is a geodesic in the (possibly disconnected) Cayley graph $\Gamma (G, X)$ from $u_1=p_-$ to $v_1=q_-$; in this case the path in $\Gamma (G, X)$ exists since $\dist_X (p_-, q_-)\leq k<K$. An analogous condition is required in the case that $u_2=p_+$. The general construction follows directly from Lemma~\ref{lem:similar}.

\emph{Paths $q_i$ and $o$, and the vertex $z$.} There exists a vertex $z$ in the segment $[v_1,v_2]$ of $q$, and a geodesic path $o$ from $u$ to $z$ satisfying the following. First
\begin{equation}\label{eq:bcp-o} \ell (o) = \dist_{X\cup \mc H}(u,z) \leq K .\end{equation}
Second, if $q_1$ and $q_2$ are the segments $[v_1, z]$ and $[z, v_2]$ of $q$ respectively, then each  $\mc H$-component of $o$ is an isolated component in one of the cycles
\begin{equation}\label{eq:cycle-1} c_1=oq_1^{-1}o_1^{-1}p_1, \ \ \ \ \  c_2 = oq_2o_2^{-1}p_2^{-1}.\end{equation}
The proof of this statement in~\cite[Proof of Prop. 3.15]{Os06} consists of a sequence of Lemmas and Corollaries:~\cite[Lem.3.16, Lem.3.17, Cor.3.18,  Lem.3.19, Cor.3.20]{Os06}. Their statements  do not involve the assumption that $X$ is a finite generating set. In the proofs, the finite generation assumption is never used.  Only~\cite[Lemma 3.16]{Os06} uses the Cayley graph $\Gamma (G, X)$ to conclude  that the path $o_1$ (respectively $o_2$) has no $\mc H$--components in the case that $u_1=p_-$ and $v_1=q_-$ (respectively $u_2=p_+$ and $v_2=q_+$);  in the case that $X$ is a finite relative generating set this also follows since $o_1$ (respectively $o_2$) was chosen to be a path in the (possibly disconnected) Cayley graph $\Gamma (G, X)$.

{\bf Final estimation and conclusion.} The constant $\epsilon_1$ is defined and it is shown that $\dist_{X\cup \Omega} (u, z) < \epsilon$:

By Proposition~\ref{thm:n-gon-2},  each $\mc H$-component $s$ of $o$ satisfies
\begin{equation} \label{eq:bcp-X} \dist_{X\cup \Omega}(s_-, s_+) \leq \dist_{\Omega} (s_-, s_+) \leq M \cdot \max\{\ell (c_1), \ell (c_2)\}. \end{equation}
By~\eqref{eq:bcp-o} and~\eqref{eq:bcp-X}, 
\begin{equation} \label{eq:bcp1-final}
\begin{split}
\dist_{X\cup \Omega}(u, z) & \leq    \ell (o)+ \ell (o) \cdot M \cdot \max\{\ell (c_1), \ell (c_2)\} \\ 
										&  \leq  \left (1+ M \cdot \max\{\ell (c_1), \ell (c_2)\} \right) \cdot K   .
\end{split}
\end{equation}
Since $q$ is a $(\mu, c)$-quasigeodesic,~\eqref{eq:bcp-p}, ~\eqref{eq:bcp-oi} and ~\eqref{eq:bcp-o} imply
\begin{equation} \label{eq:bcp-q}
\begin{split}
\ell (q_i) & \leq \mu \cdot \dist_{X\cup \mc H}(v_i, z)+c \leq \mu \cdot \left ( \ell (o_i) +\ell (p_i) +\ell (o) \right ) +c \\
		 &  \leq 2 \mu^2 K + 2\mu K + c\mu +c, \  \  \  \  \    i=1,2.
\end{split}
\end{equation}
By~\eqref{eq:cycle-1}, and~\eqref{eq:bcp-q}, 
\begin{equation}\label{eq:bcp-cycle-bound}  \ell (c_i) \leq 2 \mu^2 K + 4 \mu K  + c \mu   + 2K  + 2 c  , \ \  i=1,2.\end{equation}
Therefore~\eqref{eq:bcp1-final} and~\eqref{eq:bcp-cycle-bound} imply 
\begin{equation}\label{eq:bcp-end}
 \dist_{X\cup \Omega}(u, z)   \leq    2 \mu^2 K^2 M + 4 \mu K^2 M  + c \mu K M   + 2K^2 M  + 2 c K M +  K  .
\end{equation}
The proof concludes by defining $\epsilon_1$ the constant on the right hand side of~\eqref{eq:bcp-end}.
\end{proof}

\begin{figure}
\begin{center} 
\includegraphics[width=0.9\linewidth]{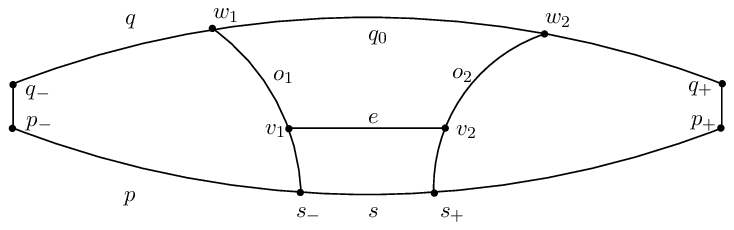}%
\end{center} 
\caption{Paths in the proof of Proposition~\ref{prop:BCP-2}.} 
 \label{fig:Figura-2} 
\end{figure}

\begin{proposition}[Bounded Coset Penetration Property, Part II]\label{prop:BCP-2}
For any $\mu \geq 1$, $c\geq 0$, $k\geq 0$, there is $\epsilon_2=\epsilon_2(\mu, c, k)$  satisfying the following condition.
Let $p$ and $q$ be a pair of $k$-similar $(\mu, c)$-quasi-geodesics in $\Gamma (G, X\cup \mc H)$ such that $p$ is without backtracking. 
 Suppose that $s$ is an $H_\lambda$--component of $p$ such that $\dist_{X\cup \Omega_\lambda}(s_-, s_+) \geq \epsilon_2$. Then
there exists an $H_\lambda$-component $t$ of $q$ such that $t$ is connected to $s$.
\end{proposition}
The proof of~\cite[Lemma 3.21]{Os06} \emph{does not} directly yield a proof of Proposition~\ref{prop:BCP-2} (it uses that the conclusion of~\cite[Proposition 3.15]{Os06} which is in terms of $\dist_X$). The argument requires an extra computation that we explain below.
\begin{proof}[Proof of Proposition~\ref{prop:BCP-2}]
Let $s$ be an $H_\lambda$--component of $p$, and suppose that there is no $H_\lambda$--component of $q$ connected to $s$. An illustration of the following construction is in Figure~\ref{fig:Figura-2}

By Proposition~\ref{prop:BCP-1}, there exists vertices $w_1, w_2$ of $q$ such that 
\begin{equation}\label{eq:sw} \{ \dist_{X\cup \Omega}(s_-, w_1) ,  \dist_{X\cup \Omega}(s_+, w_2) \} \leq \epsilon .\end{equation}

Let $r_1$ and $r_2$ paths in $\Gamma (G, X\cup \Omega)$ with $(r_1)_-=s_-$, $(r_1)_+=w_1$, $(r_2)_-=s_+$, $(r_2)_+=w_2$,
and such that
\begin{equation}\label{bcp2-r} \ell (r_i) \leq \epsilon, \  \  i=1,2.\end{equation}

Define vertices $v_1$  and $v_2$ of $r_1$ and $r_2$ as follows. 
If $r_1$ (respectively $r_2$) has an $H_\lambda$-component $t$ with $t_-=s_-$ (respectively $t_-=s_+$)
then let $v_1=t_+$ (respectively $v_2=t_+$); otherwise, let $v_1=s_-$ (respectively $v_2=s_+$). 

For $i=1,2$, let $o_i$ be the segment $[v_i, w_i]$ of the path $r_i$. By~\eqref{bcp2-r},
\begin{equation}\label{bcp2-o} \ell (o_i) \leq \epsilon, \  \  i=1,2.\end{equation}

Let $e$ be an edge from $v_1$ to $v_2$ labelled by an element of $H_\lambda$, 
and let $q_0$ be the segment $[w_1, w_2]$ of the path $q$.  Since $q$ is a $(\mu, c)$-quasi-geodesic, \eqref{bcp2-r} implies
\begin{equation} \label{bcp2-q0}
\begin{split}
 \ell (q_0) & \leq  \mu \cdot \dist_{X\cup \mc H}(w_1, w_2) + c \\
 & \leq \mu \cdot \left ( \ell (r_1) + \ell (s) +\ell (r_2) \right )  + c\\
 & \leq 2\mu \epsilon +\mu + c.
\end{split}
\end{equation}

Since $s_-$ and $v_1$ are vertices of the path $r_1 \subset \Gamma (G, X\cup \Omega)$, and
analogously $s_+$  and $v_2$ are vertices of $r_2 \subset \Gamma (G, X\cup \Omega)$,  inequality \eqref{bcp2-r} implies 
\begin{equation}\label{bcp2-v} \max\{  \dist_{X\cup \Omega}(s_-, v_1), \dist_{X\cup \Omega}(s_+, v_2) \}  \leq \epsilon. \end{equation}
From inequality~\eqref{bcp2-v},
\begin{equation} \label{bcp2-final}
\begin{split} 
 \dist_{X\cup \Omega}(s_-, s_+) & \leq \dist_{X\cup \Omega}(s_-, v_1) + \dist_{X\cup \Omega}(v_1, v_2) + \dist_{X\cup \Omega}(s_+, v_2) \\
											   & \leq 2\epsilon +\dist_{X\cup \Omega}(v_1, v_2).
\end{split}
\end{equation}

Consider the cycle
\begin{equation}\label{bcp2-cycle}  c=eo_2q_0^{-1}o_1^{-1}.\end{equation}
Observe that  $e$ is an $H_\lambda$--component of $c$. Indeed, by the choice of $o_i$, there are no $H_\lambda$-components of $o_i$ connected to $e$.
By assumption, $q$ has no $H$-components connected to $s$ (and hence $e$), and therefore $q_0$ has no $H_\lambda$-components connected to $e$. 
By Proposition~\ref{thm:n-gon-2},  
\begin{equation} \label{bcp2-ineq}
\begin{split}
\dist_{X\cup \Omega}(v_1, v_2) & \leq  \ell (c) \cdot M  \\
	& \leq \left (\ell (e) + \ell (o_2)+\ell (q_0)+\ell (o_1) \right ) \cdot  M \\
	& \leq  \left (   1 + 2 \epsilon + 2\mu \epsilon +\mu + c \right ) \cdot M,
\end{split}
\end{equation}
where the second inequality follows from~\eqref{bcp2-cycle}, and the third one from~\eqref{bcp2-o} and~\eqref{bcp2-q0}. 
Finally, ~\eqref{bcp2-final}  and~\eqref{bcp2-ineq} imply
\begin{equation}\label{eq:bcp2-end} \dist_{X\cup \Omega}(s_-, s_+)  \leq  2\epsilon +  \left (   1 + 2 \epsilon + 2\mu \epsilon +\mu + c \right ) \cdot M. \end{equation}
The proof concludes by defining  $\epsilon_2$ as the right hand side of~\eqref{eq:bcp2-end}.
\end{proof}

\begin{proposition}[Bounded Coset Penetration Property, Part III]\label{prop:BCP-3}
For any $\mu \geq 1$, $c\geq 0$, $k\geq 0$, there is $\epsilon_3=\epsilon_3(\mu, c, k)$ satisfying the following condition.
Let $p$ and $q$ be a pair of similar $k$-quasi-geodesics without backtracking in $\Gamma (G, X\cup \mc H)$. 
If $s$ and $t$ are connected $H_\lambda$--components of $p$ and $q$ respectively, then
\[\max\{ \dist_{X\cup \Omega_\lambda}(s_-, t_-), \dist_{X\cup \Omega_\lambda}(s_+, t_+) \} \leq \epsilon_3.\]
\end{proposition}
\begin{proof}
Replace in~\cite[Proof of Lemma 3.22]{Os06} the sentence ``Let us estimate the $X$-length of $c'\ $" by 
``Let us estimate the $(X\cup \Omega)$-length of $c'\ $", and use Proposition~\ref{prop:BCP-2} instead of~\cite[Lemma 3.21]{Os06}.
Then the proof goes through verbatim.
\end{proof}

\begin{lemma}\label{lem:aux-2}
Let $X,Y$ be finite relative generating sets of $G$ with respect to $\mc H$, and suppose that $X \subset Y$. Let
\[ \lambda = 1+\max \{ \dist_{X\cup \mc H}(1, y) : y \in Y  \}.\]
Then the natural inclusion $\Gamma (G, X\cup \mc H)  \hookrightarrow \Gamma (G, Y\cup \mc H)$
is a $(\lambda, 0)$-quasi-isometry that maps paths without backtracking onto paths without backtracking.
\end{lemma}
\begin{proof}
Observe that if $g=z_1z_2\cdots z_k$ where each $z_i \in Y\cup \bigcup \mc H$ and $k=\dist_{Y\cup \mc H}(1, g)$, then
\begin{equation*}  
\dist_{Y\cup \mc H}(1, g) \leq  \dist_{X\cup \mc H}(1, g) \leq \sum_{i=1}^k \dist_{X\cup \mc H}(1, z_i)  \leq \lambda \cdot \dist_{Y\cup \mc H}(1, g).
\end{equation*}
Therefore the inclusion $\Gamma (G, X\cup \mc H)\hookrightarrow \Gamma (G, Y\cup \mc H)$ is a quasi-isometric embedding,
and since the image contains all vertices of $\Gamma (G, Y\cup \mc H)$, it is in fact a quasi-isometry. The second statement on backtracking is immediate. 
\end{proof}

\subsection{Countable Relatively Hyperbolic Groups}

In this section we assume that $G$ is countable. Since countable groups can be embedded in finitely generated groups, $G$ always admits a left invariant metric $\dist:G\times G \rightarrow  \mathbb{R}$ which is proper in the sense that balls are finite. In this subsection, we state a version of the Bounded Coset Penetration property in terms of left invariant metrics as Theorem~\ref{thm:BCP-countable}. 

\begin{lemma}\label{lem:proper-left}
Let $\dist_1$ and $\dist_2$ be proper left invariant metrics on $G$. Then there is a function
$\rho : \N \rightarrow \N$ such that  $ \dist_1 (u, v) \leq \rho \left ( \dist_2 (u, v) \right ) $.
\end{lemma}
\begin{proof}
Let  $ \rho (n) = \sup \{\dist_1 (u, v) : u\in G, \  v \in G,\   \dist_2 (u, v) \leq n \}$, and observe that $\rho (n) < \infty$ for every $n \in \N$ since $\dist_2$ is proper and left invariant.
\end{proof}

\begin{lemma}\label{lem:aux-1}
Let $\dist$ be a proper left invariant metric on $G$, and $Z$ a finite subset of $G$. 
Then there is a constant  $A>0$ depending on $Z$ and $\dist$ with the following property. For any pair of elements $u$ and $v$ of $G$, if $\dist_{Z}(u, v) < \infty$, then $\dist (u,v)  \leq  A\cdot \dist_{Z}(u, v)$.
\end{lemma}
\begin{proof}
If $Z=\emptyset$ then let $A=0$. Otherwise, let
$ A = \max \{ \dist (1, z) : z \in Z \}$, and observe that $A<\infty$  since $Z$ is a finite set. 
If $\dist_Z(1, g) < \infty$, then $g=z_1\cdots z_k$ where each $z_i \in Z$
and $k=\dist_Z(1, g)$. It follows that 
\[ \dist (1, g) \leq \sum_{i=1}^k \dist (1, y_i) \leq A \cdot  \dist_{Z}(1, g) . \qedhere \]
\end{proof}

A more general version of Definition~\ref{def:similar}:
\begin{definition}[$k$-similar]
Let $\dist$ be a proper left invariant metric on $G$. For $k>0$, two paths $p$ and $q$ in $\Gamma (G, X\cup \mc H)$ are \emph{$k$-similar with respect to $\dist$} if $\max\{ \dist (p_-, q_-),  \dist (p_+, q_+) \}  \leq k$.
\end{definition}

\begin{theorem}[Bounded Coset Penetration Property]\label{thm:BCP-countable}
Let $\dist$ be a proper left invariant metric on $G$. For $\mu \geq 1, c \geq 0, k\geq 0$, there exists $\varepsilon=\varepsilon (\mu, c, k, \dist) > 0$ with the following property.   Let $p$ and $q$ be $(\mu, c)$-quasi-geodesics without backtracking in $\Gamma (G, X\cup \mc H)$.
If $p$ and $q$ are $k$-similar with respect to $\dist$, then the following conditions hold.
\begin{enumerate}
\item The sets of phase vertices of the paths $p$ and $q$ are contained in the closed $\varepsilon$-neighborhoods with respect to the metric $\dist$ of each other.
\item Let  $H\in \mc H$. If $s$ is an $H$--component of $p$ such that $\dist (s_-, s_+) \geq \varepsilon $, then there exists an $H$--component  $t$ of $q$ which is connected to $s$.
\item Let $H\in \mc H$. If $s$ and $t$ are connected $H$--components of $p$ and $q$ respectively, then 
$\max\{ \dist (s_-, t_-), \dist (s_+, t_+) \} \leq \varepsilon.$
\end{enumerate}
\end{theorem}
\begin{proof}
Let $Y=X\cup \{ g \in G : \dist (1, g) \leq k\}$. By Lemma~\ref{lem:aux-2}, there are constants $\mu'\geq 1$  and $c'\geq 0$ with the following property.
If $p$ and $q$ are $(\mu, c)$-quasi-geodesics without backtracking in $\Gamma (G, X\cup \mc H)$ which are $k$-similar with respect to $\dist$, then $p,q$ are $(\mu', c')$ -quasi-geodesics without backtracking in $\Gamma (G, Y\cup \mc H)$ which are $1$-similar with respect to $\dist_{Y}$.

Let $\Omega$ be the finite subset of $G$ provided by Proposition~\ref{thm:n-gon-2} for $G$, $\mc H$, and the finite relative generating set $Y$. 
By Lemma~\ref{lem:aux-1}, there is a constant $A>0$ such that for each $g,f \in G$
\[ \dist (g, f) \leq A \cdot \dist_{Y\cup \Omega} (g, f)  \]
whenever $\dist_{Y\cup \Omega} (g, f) < \infty$.

The conclusions of the Theorem follow from Theorem~\ref{thm:BCP} applied to $\Gamma (G, Y\cup \mc H)$ and defining
\begin{equation*} \varepsilon (\mu, c, k, \dist) =  \max \left \{ A \cdot \epsilon ( \mu', c', 1) , \   \frac{ \epsilon (\mu', c', 1) }{ A}  \right \}.  \qedhere  \end{equation*}
\end{proof}

\section{Relatively Quasiconvex Subgroups} \label{sec.quasiconvexity}

Let $G$ be a countable group hyperbolic relative to a finite collection of subgroups $\mc H$, and let  $X$  be a finite relative generating set of $G$ with respect to $\mc H$.

\subsection{Definition of Quasiconvex Subgroups}

\begin{definition}[Parabolic Subgroups]
A subgroup $P$ of $G$ is called \emph{parabolic with respect to $\mc H$} if there is $H \in \mc H$ and $g \in G$ such that $P<gH  g^{-1}$. 
\end{definition}

\begin{definition}[Relatively Quasiconvex Subgroup]\label{def:qc}
A subgroup $Q$ of $G$ is called \emph{quasiconvex relative to $\mc H$} or simply \emph{quasiconvex in $(G, \mc H)$}  if 
there exists a constant $\sigma \geq 0$ such that the following holds: Let $f$, $g$ be two elements of $Q$, and $p$ an arbitrary geodesic path from $f$ to $g$ in  $\Gamma(G, X\cup \mc H)$.  For any (phase) vertex $v \in p$, there exists a vertex $w \in Q$ such that $ \dist (v,w) \leq \sigma$.

The definition relative quasiconvexity is independent of the finite relative generating set $X$, and the proper left invariant metric $\dist$, see Proposition~\ref{prop:qc-independence} below. In particular, we say \emph{$Q$ is a quasiconvex subgroup of $(G, \mc H)$} to express that $Q$ is a  quasiconvex subgroup of $G$ relative to $\mc H$. If referring to the constant $\sigma$ is necessary, we say that \emph{$Q$ is $\sigma$-quasiconvex in $(G, \mc H, X, \dist)$}.
\end{definition}

\begin{remark}
Definition~\ref{def:qc} appeared in~\cite{Os06} assuming that $G$ is finitely generated, and the version for countable groups appeared in~\cite{HK08}.  Definition~\ref{def:qc} is equivalent to other approaches to relative quasiconvexity and we refer the interested reader to~\cite{HK08, MaWi10b} and the references there in.
\end{remark}

That Definition~\ref{def:qc} of relative quasiconvexity is independent of $X$ and $\dist$, and that relative quasiconvexity is preserved by conjugation are proved in~\cite{HK08} passing through an equivalent definition of relative quasiconvexity. We provide direct proofs using Theorem~\ref{thm:BCP-countable}.

\begin{proposition}\cite{HK08} \label{prop:qc-independence}
Definition~\ref{def:qc} is independent of the choice of finite relative generating set $X$ and the choice of proper metric $\dist$. 
\end{proposition}
\begin{proof}
Let $Q$ be a subgroup of $G$,  
let $\dist$ be a proper left invariant metric on $G$,  let $X$ and $Y$ be finite relative generating sets of $G$ with respect to $\mc H$,
and suppose that $X\subset Y$. By Lemma~\ref{lem:aux-2} and Theorem~\ref{thm:BCP-countable}, there is $\tau >0$ 
such that if $p$ is a geodesic in  $\Gamma (G, Y\cup \mc H)$ and $q$ is a geodesic in $\Gamma (G, X\cup \mc H)$ between the same points of $Q$, 
then the sets of  vertices of $p$ and $q$ are contained in the closed $\tau$-neighborhoods  with respect to $\dist$ of each other.
Hence $Q$ is relatively quasiconvex with respect to  $X$ and $\dist$ if and only if $Q$ is relatively quasiconvex with respect to $Y$
and $\dist$. Since the union of two finite relative generating sets  is a finite relative generating set,  Definition~\ref{def:qc} is independent of the choice of finite relative generating set.

Let $X$ be a finite relative generating set of $G$ with respect to $\mc H$, and let $\dist_1$ and $\dist_2$ be proper left invariant metrics on $G$. 
By Lemma~\ref{lem:proper-left}, there is a function $\rho : \N \rightarrow \N$ such that 
$ \dist_1 (u, v) \leq \rho \left ( \dist_2 (u, v) \right ) $ for any $u$ and $v$ in $G$. Therefore, if $Q$ is relatively quasiconvex with respect to the relative generating set $X$ and $\dist_2$, then $Q$ is relatively quasiconvex with respect to $X$ and $\dist_1$. Since $\dist_1$ and $\dist_2$ were arbitrary, Definition~\ref{def:qc} is independent of the choice of proper left invariant metric.
\end{proof}

\begin{corollary}\cite{HK08} \label{prop:qc-conj}
If a subgroup $Q<G$ is quasiconvex relative to $\mc H$, then any conjugate of $Q$ is quasiconvex relative to $\mc H$.
\end{corollary}
\begin{proof}
Suppose that $Q$ is $\sigma$--quasiconvex in $\Gamma (G, X\cup \mc H)$,  and let $g \in G$. By Proposition~\ref{prop:qc-independence}, without loss of generality, we can assume that $g \in X$. Observe that since $Q$ is $\sigma$--quasiconvex, if $q$ is an arbitrary geodesic with endpoints in the left coset $gQ$, then for any (phase) vertex $v \in q$, there exists a vertex $w \in gQ$ such that $ \dist (v,w) \leq \sigma$. 

Let $p$ be a geodesic in $\Gamma (G, X\cup \mc H)$ with endpoints in $gQg^{-1}$, and let $u$ be a (phase) vertex on $p$.
Let $q$ be a geodesic between the vertices $p_-g$ and $p_+g$.  Since $g \in X$, the geodesic paths $p$ and $q$ are $1$-similar with respect to $\dist$. 
Therefore Theorem~\ref{thm:BCP-countable} implies that there is a (phase) vertex $v$ on $q$ such that $\dist (u, v) \leq \epsilon (1,0,1)$.
Since $q$ is a geodesic with endpoints in $gQ$, there is a vertex $w \in gQ$ such that $\dist (v, w) \leq \lambda$. It follows that
the $wg^{-1} \in gQg^{-1}$ and
\begin{equation*}
\dist (u, wg^{-1}) \leq \dist (u, v) + \dist (v, w) + \dist (w, wg^{-1})  \leq \epsilon (1,0,1) + \lambda + 1. 
\end{equation*}
This shows that $gQg^{-1}$ is quasiconvex in $\Gamma (G, X\cup \mc H)$.
\end{proof}

\begin{corollary}\cite{HK08}\label{cor:par-qc}
A parabolic subgroup of $(G, \mc H)$ is quasiconvex in $(G, \mc H)$.
\end{corollary}
\begin{proof}
Let $Q$ be a parabolic subgroup. Then there is $H \in \mc H$, $P <H$, and $g \in G$ such that $Q=gPg^{-1}$. Observe that any geodesic in $\Gamma (G, X\cup \mc H)$ with endpoints in $P$ is a single edge, and therefore $P$ is quasiconvex relative to $\mc H$. 
By Corollary~\ref{prop:qc-conj}, $Q$ is quasiconvex relative to $\mc H$.
\end{proof}

\subsection{Properties of Relatively Quasiconvex Subgroups} 

\begin{theorem} \cite{HK08,MP09} \label{prop:intersection}
Let $Q$ and $R$ be quasiconvex subgroups of $G$ relative to $\mc H$. Then $Q\cap R$ is quasiconvex in $G$ relative to $\mc H$.
\end{theorem}

\begin{theorem} \cite{HK08, MaWi10b}  \label{thm:qc-rh}
Let $Q$ be a $\sigma$-quasiconvex subgroup of $G$ relative to $\mc H$, and let $\mc L = \{ Q\cap gHg^{-1}   :  g\in G, H\in \mc H, \dist (1, g) \leq \sigma\}$.
Then $Q$ is hyperbolic relative to $\mc L$.  In particular, $Q$ is finitely generated relative to $\mc L$, and each infinite maximal parabolic subgroup of $Q$ with respect to $\mathcal{H}$ is conjugate in $Q$ to an element of $\mc L$
\end{theorem}

\section{Quasigeodesics and Hyperbolically embedded subgroups} \label{sec.quasigeodesics}

Let $G$ be a countable group hyperbolic relative to a finite collection of subgroups $\mc H$, let  $X$  be a finite relative generating set of $G$ with respect to $\mc H$,
let $\dist$ be a proper left invariant metric on $G$. 

Hyperbolically embedded subgroups were defined in the introduction, see Definition~\ref{def.hyp_embedded}.  If $P$ is a hyperbolically embedded subgroup in $(G, \mc H)$ then we denote by $\Gamma (G, X\cup \mc H \cup \{P\})$ the Cayley graph of $G$ with respect to the generating set $P \sqcup \bigsqcup_\Lambda H_\lambda$.

\subsection{Hyperbolically Embedded Subgroups}

The hyperbolically embedded subgroups where characterized by Bowditch~\cite[Theorem 7.11]{BO99} for the case that $G$ is a hyperbolic group, and by Osin~\cite[Theorem 1.5]{Os06-3} for the general case:
\begin{theorem}\cite[Theorem 1.5]{Os06-3}\label{thm:hyp-emb-char}
A subgroup $P$ of $G$ is hyperbolically embedded into $(G, \mc H)$, if and only if the following conditions hold. 
\begin{enumerate}
 \item\label{fin-gen} $P$ is generated by a finite set $Y$.
 \item\label{quasiconvex} There exists $\lambda, c \geq 0$ such that for any element $g\in P$, we have $|g|_Y \leq \lambda \cdot \dist_{X\cup \mc H}(1, g) + c$.
 \item\label{malnormal} For any $g \in G$ such that $g \not \in P$, the intersection $P\cap P^g$ is finite.
\end{enumerate}
\end{theorem}

\begin{corollary}\cite[Theorem 1.5]{Os06-3} \label{thm:hyp-qcx}
Let $P$ be a hyperbolically embedded into $(G, \mc H)$. Then $P$ is a hyperbolic group, and $P$ is quasiconvex subgroup of $(G, \mc H)$.
\end{corollary}
\begin{proof}
By Theorem~\ref{thm:hyp-emb-char}, $P$  has a finite generating set $Y$. Without loss of generality, assume that $Y \subset X$.
Then the Cayley graph of $\Gamma (P, Y)$ is quasi-isometrically embedded into the hyperbolic space $\Gamma(G, X\cup \mc H)$, and therefore $P$ is a hyperbolic group.  Since geodesics in $\Gamma (P, Y)$ are $(\lambda, c)$-quasi-geodesics in $\Gamma(G, X\cup \mc H)$, Theorem~\ref{thm:BCP-countable} implies that vertices of geodesics in $\Gamma(G, X\cup \mc H)$ with endpoints in $P$ stay closed to elements of $P$ with respect to $\dist$. Therefore $P$ is quasiconvex in $(G, \mc H)$.
\end{proof}

The proof of Corollary~\ref{cor:sandwich} is based on the argument due to Osin in the proof of the ``only if" part of Theorem~\ref{thm:hyp-emb-char} in~\cite{Os06-3}. 
\begin{corollary}\label{cor:sandwich}
Let $(G, \mc H_1)$ and $(G, \mc H_2)$ be relatively hyperbolic structures with $\mc H_1 \subset \mc H_2$ and $\mc H_2$ finite. Then for each collection of subgroups $\mc H$ such that $\mc H_1 \subset \mc H \subset \mc H_2$, the pair $(G, \mc H)$ is relatively hyperbolic.  In particular, each subgroup $P$ in $\mc H_2 \setminus \mc H_1$ is hyperbolically embedded into $(G, \mc H_1)$. 
\end{corollary}
\begin{proof}
It is enough to show that every element of $\mc H_2 \setminus \mc H_1$ is  hyperbolically embedded into $(G, \mc H_1)$. Then the statement of the corollary follows by induction on the size of $\mc H \setminus \mc H_1$.  

Suppose that $X$ is a finite relative generating set of $(G, \mc H_1)$. For each $P \in \mc H_2 \setminus \mc H_1$, let $\Omega_P$ be the subset given by Proposition~\ref{thm:n-gon-2} for the relatively hyperbolic pair $(G, \mc H_2)$.  Let $M$ be the constant of Proposition~\ref{thm:n-gon-2} for the pair $(G, \mc H_2)$.  Then for each element $g \in P$, one considers a geodesic path $q$ in $\Gamma (G, X\cup \mc H_1)$ from $1$ to $g$.  Using that $\Gamma (G, X\cup \mc H_1) \subset \Gamma (G, X\cup \mc H_2)$ and Proposition~\ref{thm:n-gon-2}, it follows that $P$ is generated by the finite subset $\Omega_P$, and that $|g|_{\Omega_P} \leq M \cdot \dist_{X\cup \mc H_1}(1, g)$. Moreover, by Theorem~\ref{thm:malnormality}, for any $f \in G \setminus P$, the intersection $P\cap P^f$ is finite. By Theorem~\ref{thm:hyp-emb-char}, $P$ is hyperbolically embedded into $(G, \mc H_1)$.
\end{proof}

\subsection{Quasigeodesics and Hyperbolically Embedded Groups}

Let $P$ be a hyperbolically embedded subgroup of $(G, \mc H)$, and let $\dist_{X\cup \mc H}$ denote the combinatorial metric of the graph $\GammaH$. 

\begin{lemma} \label{lem:hyp_emb}
There exist a constant $\mu >0$ with the following property. Let $p$ be a path in $\Gamma (G, X\cup \mc H\cup \{P\})$ without backtracking. Let $q$ be a path in $\Gamma (G, X\cup \mc H)$ obtained by replacing each $P$-component of $p$ by a geodesic segment in $\Gamma (G, X\cup \mc H)$ between its endpoints.  Then $\displaystyle  \ell (q) \leq \mu \cdot \left( \ell (p) + \dist_{X\cup \mc H}(q_-, q_+) \right)$. 

In particular, if $p$ is a $(\lambda, c)$-quasigesodesic in $\Gamma (G, X\cup \mc H\cup \{P\})$,  then \[ \ell (q) \leq \mu \cdot (\lambda+1) \cdot  \dist_{X\cup \mc H} (q_-, q_+) + \mu c. \]
\end{lemma}
\begin{proof}
Since $P$ is hyperbolically embedded, the group $G$ is hyperbolic relative to $\mc H\cup \{P\}$.
Let $M>0$ and $\Omega \subset G$ be given by Proposition~\ref{thm:n-gon-2} for the pair $(G, \mc H\cup \{P\})$ and the relative generating set $X$. 

Suppose that $p$ is of the form $p = s_1 t_1 \dots s_k t_k s_{k+1}$ where $t_1, \dots , t_k$ are all the $P$-components of $p$.
Then $q = s_1 u_1 \dots s_k u_k s_{k+1}$ where each $u_i$ is a geodesic segment in $\GammaH$ connecting the endpoints of $t_i$. 
First observe that 
\begin{equation*} 
 \ell (q) \leq  \ell (p) +  \sum_{i=1}^{k} \ell (u_i). 
\end{equation*}

Let $r$ be a geodesic in $\GammaH$ from $p_+$ to $p_-$.  Consider the closed polygon in $\Gamma (G, X\cup \mc H\cup \{P\})$ given by $p r$.
Specifically this cycle decompose as $p r = s_1 t_1 \dots s_k t_k s_{k+1} r$.  Since $p$ has not backtracking,  $t_1, \dots , t_k$ is a collection of isolated $P$-components of this cycle. Therefore, if $g_1, \dots , g_k$ are the elements of $G$ represented by the labels of $t_1, \dots , t_k$, then Proposition~\ref{thm:n-gon-2} implies 
\begin{equation*} 
  \sum_{i=1}^{k}  |g_i|_{X\cup \Omega}  \leq  M \cdot ( \ell (p)+\ell (r) ) 
\end{equation*}
Since $\ell (u_i)=|g_i|_{X \cup \mc H} \leq |g_i|_{X\cup \Omega}$ for each $g_i$, and $\ell (r) = \dist_{X\cup \mc H}(q_-, q_+)$, we have
that 
\begin{equation*} \sum_{i=1}^{k} \ell (u_i)  \leq M \cdot \left( \ell (p) + \dist_{X\cup \mc H}(q_-, q_+) \right). \end{equation*} 
Therefore $\ell (q) \leq 2M \cdot \left ( \ell (p) + \dist_{X\cup \mc H}(q_-, q_+) \right )  $ which is the first conclusion of the lemma. 

The second statement follows by observing that the combinatorial metric in $\GammaH$ is always larger than the one in $\Gamma (G, X\cup \mc H\cup \{P\})$,
and that $p$ and $q$ have the same endpoints.
\end{proof}

\begin{proposition} \label{prof:hyp_emb}
For $\lambda \geq 1$ and $c \geq 0$, there is a constant $L = L (\lambda, c) \geq 1$ with the following property. If $p$ is a $(\lambda, c)$-quasi-geodesic in $\Gamma (G, X\cup \mc H \cup \{P \})$ and $q$ is a path in $\GammaH$ obtained by replacing each $P$-component of $p$ by a geodesic segment in $\GammaH$ connecting its endpoints, then $q$ is a $(L, L)$-quasi-geodesic in $\GammaH$.
\end{proposition}
\begin{proof}
Let $\mu$ be the constant given by Lemma~\ref{lem:hyp_emb} and let $L > 4\mu (1+\lambda + c)$.
Suppose that $p = s_1 t_1 \dots s_k t_k s_{k+1}$ where $t_1, \dots , t_k$ are all the $P$-components of $p$. Then $q = s_1 u_1 \dots s_k u_k s_{k+1}$ where each $u_i$ is a geodesic segment in $\GammaH$ connecting the endpoints of $t_i$.  Let $q'$ be a subpath of $q$. Below we show that $\ell (q') \leq L \cdot \dist_{X\cup \mc H}(q'_-, q'_+) + L$ in two cases: 

{\it Case 1.} Suppose that $q'$ is of the form $ q' = s'_i u_i \dots s_j u_j s'_{j+ 1}$, where $s'_i$ and $s'_{j + 1}$ are subpaths of $s_i$ and $s_{j + 1}$ respectively.  In this case, Lemma ~\ref{lem:hyp_emb} implies that
\[ \ell (q') \leq \mu  \cdot (\lambda+1) \cdot  \dist_{X\cup \mc H}(q'_-, q'_+) + \mu c. \]

{\it Case 2.} Otherwise, at least one of the endpoints of $q'$ is a vertex of $u_i$ for some $i$.
This yields three similar subcases. We only consider one and leave the other two for the reader.
Suppose that $q' = u'_i s_{i+1} u_{i +1} \dots u_{j-1} s_{j} u'_{j}$ where $u'_i$ and $u'_{j}$ are subpaths of $u_i$ and $u_{j}$ respectively.
Let $p'$ be the path in $\Gamma (G, X\cup \mc H \cup \{P \} )$ given by $p' = t'_i s_{i+1} t_{i +1} \dots t_{j-1} s_{j} t'_{j}$,
where $t'_i$ and $t'_{j}$ correspond to single edges connecting the endpoints of $u'_i$ and $u'_{j}$ respectively. 
Since the subpath $s_{i+1} t_{i +1} \dots t_{j-1} s_{j}$ of $p'$ is a $(\lambda ,c)$-quasi-geodesic without backtracking, it follows that 
$p'$ is a $(\lambda, 4\lambda+c)$-quasi-geodesic without backtracking. Lemma ~\ref{lem:hyp_emb} implies that
\[ \ell (q') \leq \mu \cdot (\lambda+1) \cdot  \dist_{X\cup \mc H}(q'_-, q'_+) + 4 \mu \lambda + \mu c .  \qedhere \]
\end{proof}

\section{Proofs of the Main Results}\label{sec:proofs} 

\subsection{Proof of Theorem~\ref{thm:main}~\eqref{main-1}} 

\begin{proposition} \label{thm:main-1}
Let $G$ be a countable group, let $(G, \mc H)$ be a relatively hyperbolic structure.  Suppose that $P$ is hyperbolically embedded subgroup of $(G, \mc H)$   and $Q$ is a quasiconvex subgroup of $(G, \mc H)$. Then $Q$ is quasiconvex subgroup of $(G, \mathcal{H} \cup \{ P \})$.
\end{proposition}
\begin{proof}
Let $X$ be a finite relative generating set of $(G, \mc H)$, and let $\dist$ be a proper left-invariant metric on $G$. Suppose that $Q$ is a $\sigma$-quasiconvex subgroup of $(G, \mc H)$ with respect to $X$ and $\dist$. Let $f$ be an element of $Q$ and let $p$ be a geodesic in $\Gamma (G, X \cup \mc H \cup \{ P\})$ from $1$ to $f$. Let $q$ be the path in $\GammaH$ obtained by replacing each $P$-component of $p$ by a geodesic segment in $\GammaH$. Let $r$ be a geodesic in $\GammaH$ from $1$ to $f$. By Proposition ~\ref{prof:hyp_emb} implies that $q$ is a $(\lambda,c)$-quasi-geodesic in $\GammaH$, where the constants $\lambda$ and $c$ are independent of the element $f$ and the path $p$.  By Theorem~\ref{thm:BCP-countable}, there is constant $\epsilon=\epsilon (\lambda, c, 0)$  such that all phase vertices of $q$ are in the $\epsilon$-neighborhood of the vertices of $r$ with respect to $\dist$.  By hypothesis  all vertices of $r$ are in the $\sigma$-neighborhood of $Q$ with respect to $\dist$, it follows that all vertices of $p$ are in the $(\epsilon+\sigma)$-neighborhood of $Q$. Therefore $Q$ is a $(\epsilon+\sigma)$-quasiconvex subgroup of $(G, \mc H \cup \{P\})$ with respect to $X$ and $\dist$.
\end{proof}

\begin{proof}[Proof of Theorem~\ref{thm:main}~\eqref{main-1}]
By Corollary~\ref{cor:sandwich} every element of $\mc{H}_2\setminus \mc{H}_1$ is hyperbolically embedded into  $(G , \mc{H}_1)$.
If $Q<G$ is quasiconvex in $(G, \mc H_1)$, an induction argument using Proposition~\ref{thm:main-1} shows that $Q$ is quasiconvex in $(G, \mc H_2)$.
\end{proof}

\subsection{Proof of Theorem~\ref{thm:main}~\eqref{main-2}}  

Let $G$ be a countable group  and let $\dist$ be a proper left invariant metric on $G$.  

\begin{convention}
For $A \subset G$ and $\alpha >0$, we denote by $N_\alpha (A)$ the closed $\alpha$-neighborhood of $A$ in the metric space $(G, \dist )$.
\end{convention}

\begin{lemma}\cite{HK08, MP09} \label{lem:quasiorthogonality1}
Let $P$ and $Q$ be subgroups of $G$. For each $g\in G$ and $\sigma \geq 0$ there is $L = L(P, Q, g, \sigma) \geq 0$ so that 
\[ Q \cap N_{\sigma}(gP)  \subset  N_L ( Q \cap gPg^{-1}).\] 
\end{lemma}
\begin{proof}
Suppose the statement is false. Then there are sequences $\{ q_n \}_{n=1}^\infty$ and $\{ h_n \}_{n=1}^\infty$ such that $q_n \in Q$, $q_n h_n \in gP$, $\dist (1,h_n) \leq \sigma$, and $\dist (q_n, Q\cap gPg^{-1} ) \geq n $.  
Since $\dist$ is proper, without loss of generality, assume $\{ h_n \}_{n=1}^\infty$ is a constant sequence $\{ h\}_{n=1}^\infty$.
Observe that $q_nq_m^{-1} = (q_nh) (q_m h)^{-1} \in Q\cap gPg^{-1}$, and hence $q_mh$ and $q_nh$
are in the same right coset of $Q\cap gPg^{-1}$, say $ (Q\cap gPg^{-1}) f$. It
follows that
\[ \dist (q_n, Q \cap gPg^{-1} ) \leq \dist (q_n, q_nh) + \dist (q_nh, Q \cap gPg^{-1}) \leq \sigma + \dist (1,f) \]
for any $n$, a contradiction.
\end{proof}

\begin{lemma}\label{lem:parabolic-aprox}
Let $P$ and $Q$ be subgroups of $G$. For each $\sigma >0$, there is  $L=L(P, Q, \sigma)>0$ with the following property. 
If $q \in Q$ is a product of the form $q =g p f$ where $\dist  (1, g) \leq \sigma$, $\dist  (1, f) \leq \sigma$, and $p \in P$, then  $q=a b$ where $a \in    Q \cap gPg^{-1}$,  $ b \in Q$, and $\dist (1, b) \leq L$. 
\end{lemma}
\begin{proof}
For $g \in G$, let $M_g =M(   P,  Q,  g, \sigma)$ be the constant provided by Lemma~\ref{lem:quasiorthogonality1}. Since $\dist$ is proper, $L = \max \{   M_g\ :\  g \in    G,\ \dist (1, g) \leq \sigma \}$ is a well-defined positive integer. Suppose that $q \in    Q$ is a product of the form $q =g p f$ where $\dist (1, g) \leq \sigma$, $\dist (1, f) \leq \sigma$, and $p \in    P$. Lemma~\ref{lem:quasiorthogonality1} implies that
\[ q \in     Q \cap N_{\sigma}(gP) \subset N_{L}(    Q \cap gPg^{-1} ).\]
The conclusion of the lemma follows.
\end{proof}

Let $\mc H$ be a finite collection of subgroups of $G$, let $X$ be a relative generating set of $(G, \mc H)$,  and suppose that $(G, \mc H)$ is a relatively hyperbolic structure. 

\begin{lemma}\label{lem:qc-proper-metric}
Let $R$ be a $\varsigma$-quasiconvex subgroup of $(G, \mc H, X, \dist)$. For any $\lambda \geq 1$, $c\geq 0$ and $L>0$, there is 
$\tau=\tau ( L, \lambda, c, \varsigma)>0$ with the following property.  

If $u$ is a $(\lambda, c)$-quasigeodesic without backtracking in $\Gamma (G, X\cup \mc H)$ such that its endpoints are in $N_L ( R)$,  
then all phase vertices of $u$ are in $N_\tau ( R)$. 
\end{lemma}
\begin{proof}
Observe that there is a geodesic $v$ in $\Gamma (G, X\cup \mc H)$ with endpoints in $Q$, and such that $u$ and $v$ are $L$-similar with respect to $\dist$.  
The Lemma follows immediately from Theorem~\ref{thm:BCP-countable} by taking $\tau=\sigma + \epsilon (\lambda, c, L, \dist)$.
\end{proof}

\begin{lemma}\label{lem:main-thm}
Let $P$ and $Q$ subgroups of $G$, and suppose there are constants $\sigma, \varsigma >0$  such that for any $g \in G$ with $\dist (1, g) \leq \sigma$, the subgroup $Q\cap gPg^{-1}$ is $\varsigma$-quasiconvex in $(G, \mc H, X, \dist )$. Then there is a constant $\tau=\tau (\sigma, \varsigma)>0$ with the following property. 
If $u$ is a geodesic in $\GammaH$ with endpoints in $N_\sigma (Q)$ and such that $\Label (u)$ represents an element of $P$, then all vertices of $u$ are in $N_\tau (Q)$. 
\end{lemma}
\begin{proof}
Let $L=L(Q, P, \dist, \sigma)>0$ be provided by Lemma~\ref{lem:parabolic-aprox}, and let $\tau=\tau ( L, 1, 0, \varsigma)>0$ be the constant provided by Lemma~\ref{lem:qc-proper-metric}. 

Let $u$ be geodesic segment as in the statement of the lemma. By hypothesis, there are elements $x,y  \in Q$ such that 
$ \label{eq:claim1-mt} \max \{ \dist ( x, u_-),  \dist ( y,u_+) \} \leq \sigma $. Let $g$ denote the element $x^{-1}u_-$, $p$ denote the element $u_-^{-1}u_+$, and $f$ the element $u_+^{-1}y$. Since $x^{-1}y=gpf$ is an element of $Q$,  Lemma~\ref{lem:parabolic-aprox} implies that $x^{-1}y=gpf=a b$ where $a \in Q \cap gPg^{-1}$,  $b \in Q$, and $\dist (1, b) \leq L$.

Let $x^{-1}u$ be the translation of $u$ by $x^{-1}$. Then $x^{-1}u$ is a geodesic segment with endpoints in $N_L ( Q\cap gPg^{-1})$.  By Lemma~\ref{lem:qc-proper-metric}, all vertices of $x^{-1}u$ are in $N_\tau ( Q\cap gPg^{-1} )$. Since $x \in Q$, all vertices of $u$ are in $N_\tau (Q)$. 
\end{proof}

\begin{lemma}[Quasiconvexity Criterion]\label{lem:qc-equiv}
Let $Q$ be a subgroup of $G$. Suppose that there are constants $\lambda \geq 1$, $c\geq 0$, $\sigma >0$ such that for any $f \in Q$ there is a $(\lambda,c)$-quasigeodesic $q$ in $\Gamma (G, X\cup \mc H)$ from the $1$ to $f$ with all phase vertices in $N_\sigma (Q)$. Then $Q$ is a quasiconvex in $(G, \mc H)$.
\end{lemma}
\begin{proof}
If $p$ is a geodesic in $\Gamma (G, X\cup \mc H)$ from $1$ to $f$, then Theorem~\ref{thm:BCP-countable} imply that all vertices of $p$ are within $(\sigma +\epsilon (\lambda, c, 0))$ of elements of $Q$ with respect to $\dist$. Hence $Q$ is quasiconvex in $(G, \mc H)$.
\end{proof}

\begin{proposition}\label{thm:main-2}
Suppose that $P$ is a hyperbolically embedded subgroup of $(G, \mc H)$  and $Q$ is a $\sigma$-quasiconvex subgroup of $(G, \mc H \cup \{P\}, X, \dist )$.  If there is $\varsigma>0$ such that for each $g \in G$ with $\dist (1, g) \leq \sigma$ the subgroup $Q\cap gPg^{-1}$ is $\varsigma$-quasiconvex in $(G, \mc H, X, \dist )$, then $Q$ is quasiconvex in $(G, \mc H)$.
\end{proposition}
\begin{proof}
Let $L$ be the constant provided by Proposition~\ref{prof:hyp_emb} for $P$.  Let $f$ be an element of $Q$,  $p$ a geodesic in $\Gamma (G, X\cup \mc H\cup \{P\} )$ from $1$ to $f$, and $q$ be the path in $\GammaH$ obtained by replacing each $P$-component of $q$ by a geodesic segment in $\GammaH$. Proposition~\ref{prof:hyp_emb} imply $q$ is an $(L,L)$-quasigeodesic in $\Gamma (G, X\cup \mc H)$. Below we show that there is a constant $\kappa =\kappa (\sigma, \varsigma)$ such that all phase vertices of $q$ are in $N_\kappa (Q)$. Then Lemma~\ref{lem:qc-equiv} implies that $Q$ is quasiconvex in $(G, \mc H)$.

Decompose the path $p$ as \[p = s_1 t_1 \dots s_k t_k s_{k+1},\]  where $\{ t_i \}_{i=1}^{k}$ are all the $P$-components of $p$. Then the path $q$ in $\GammaH$ decomposes as \[q = s_1 u_1 \dots s_k u_k s_{k+1},\]  where each $u_i$ is a geodesic segment in $\GammaH$ connecting the endpoints of $t_i$. 

Since $Q$ is $\sigma$-quasiconvex in $(G, \mc H \cup \{P\}, X , \dist)$, all phase vertices of $s_i$ are in $N_\sigma (Q)$ and, in particular, the endpoints of each $u_i$ are  in $N_\sigma (Q)$. Since $\Label (u_i)$ represents an element of $P$ and $Q\cap gPg^{-1}$ is $\varsigma$-quasiconvex in $(G, \mc H, X, \dist )$ for each $g \in G$ with $\dist (1, g) \leq \sigma$,  Lemma~\ref{lem:main-thm} implies that there is $\tau=\tau (\sigma, \varsigma)>0$ such that all vertices of each $u_i$ are in $N_\tau (Q)$.  It follows that all vertices of $q$ are in $N_\kappa (Q)$ where $\kappa = (\tau + \sigma)$. 
\end{proof}

\begin{proof}[ Proof of Theorem~\ref{thm:main}~\eqref{main-2} ]
Suppose that $\mc{H}_2 \setminus \mc{H}_1= \{ P_1, \dots , P_l \}$ and let $\mc P_k = \mc H_1 \cup \{P_i : 1 \leq i \leq k\}$.  Suppose that $Q$
is quasiconvex in $(G, \mc{H}_2)$ and that for each $g\in G$ and $P_i \in \mc{H}_2\setminus \mc{H}_1$, the subgroup $Q\cap g P_i g^{-1}$ is quasiconvex in  $(G, \mc \mc{H}_1)$. 

Corollary~\ref{cor:sandwich} implies that $(G, \mc P_{k-1})$ is a relatively hyperbolic structure, and that $P_k$ is hyperbolically embedded into $(G, \mc P_{k-1})$.  Since for each $g \in G$  the subgroup $Q\cap g P_k g^{-1}$ is quasiconvex in $(G, \mc H_1)$, Proposition~\ref{thm:main-1} implies that $Q\cap g P_k g^{-1}$ is   quasiconvex in  $(G, \mc P_{k-1})$.  This shows that the hypothesis of Proposition~\ref{thm:main-2} are satisfied for $Q$, $P_k$ and $(G, \mc P_{k-1})$, and therefore $Q$ is quasiconvex in $(G, \mc P_{k-1})$. An induction argument shows that $Q$ is quasiconvex in $(G, \mc{H}_1)$.
\end{proof}

\subsection{Proof of Corollary~\ref{cor:elementary}}

\begin{lemma}\label{lem:elem}
If $H$ is a hyperbolic group with the property that any subgroup is finitely generated, then $H$ is an elementary subgroup. 
\end{lemma}
\begin{proof}
Since every non-elementary hyperbolic group contains a subgroup isomorphic to a free group in two generators, the proposition is immediate. 
\end{proof}

\begin{proof}[Proof of Corollary~\ref{cor:elementary}]
Assume that $(G, \mc H_1)$ and $(G, \mc H_2)$ have the same class of quasiconvex subgroups, and let $H$ be any subgroup in $\mc H_2 \setminus \mc H_1$.
By Corollaries~\ref{thm:hyp-qcx} and~\ref{cor:sandwich}, $H$ is hyperbolically embedded into $(G, \mc H_1)$  and therefore is a hyperbolic group. 
By the assumption and Corollary~\ref{cor:par-qc},  each subgroup of $H$ is quasiconvex in $(G, \mc H_1)$.  Therefore, Theorem~\ref{thm:qc-rh} implies that every subgroup of $H$ is finitely generated relative to $\mc H_1$. By Theorem~\ref{thm:malnormality}, the intersection of $H$ with any conjugate in $G$ of a subgroup in $\mc H_1$ is finite. It follows that any subgroup of $H$ is finitely generated. By Lemma~\ref{lem:elem}, $H$ is an elementary group.

Suppose that each subgroup in $\mc H_2 \setminus \mc H_1$ is elementary. Since elementary subgroups are always quasiconvex~\cite[Corollary 1.7]{Os06-3}, Theorem~\ref{thm:main} implies that  $(G, \mc H_1)$ and $(G, \mc H_2)$ have the same class of quasiconvex subgroups.
\end{proof}

\subsection{Proof of Theorem~\ref{thm:coh} }

Suppose that each $H \in \mc H$ is coherent, and let $Q$ be a finitely generated subgroup of $G$. By hypothesis, $Q$
is quasiconvex in $(G, \mc H)$, and therefore, by Theorem~\ref{thm:qc-rh}, hyperbolic relative to a finite collection $\mc L$ of parabolic subgroups of $G$. 
It follows that each $L \in \mc L$ is finitely generated~\cite[Proposition 2.29]{Os06}, and by coherence of the maximal parabolic subgroups of $G$, each $L \in \mc L$ is
finitely presented. By~\cite[Corollary 2.41]{Os06}, $Q$ is finitely presented.

\subsection{Proof of Corollary~\ref{cor:coh} }

Suppose that each maximal parabolic subgroup of $G$ is a locally quasiconvex hyperbolic group. Since each $H\in \mc H$ is hyperbolic, the group $G$ is hyperbolic~\cite[Corollary 2.41]{Os06}. Let $Q<G$ be finitely generated. By hypothesis, $Q$ is quasiconvex in $(G, \mc H)$, and therefore, by Theorem~\ref{thm:qc-rh}, hyperbolic relative to a finite collection $\mc L$ of parabolic subgroups of $G$.  Each $L \in \mc L$ is finitely generated by~\cite[Proposition 2.29]{Os06}. By local quasiconvexity of the maximal parabolic subgroups of $G$, each $L \in \mc L$ is quasiconvex in a maximal parabolic subgroup of $G$, and hence quasiconvex in $G$. By Theorem~\ref{thm:main}, $Q$ is quasiconvex in $G$.

\subsection{Proof of Corollary~\ref{cor:3}}\label{subsec:qc-core}

We describe almost verbatim the peripheral structure induced by a quasiconvex subgroup in a hyperbolic group following the exposition in~\cite[Section 3.1]{AGM08}. 
Let $G$ be a hyperbolic group, and let $Q$ be a quasiconvex subgroup of $G$.

\begin{lemma}[Malnormal Core of $Q$] \label{lem:qc-core}
There is a finite collection $\mc D$ of infinite quasiconvex subgroups of $Q$ with the following properties.
\begin{enumerate}
\item For each $D \in \mc{D}$ and $g \in Q\setminus D$, the intersection $D \cap D^g$ is finite, and
\item For each distinct pair $D, D' \in \mc{D}$ and $g \in Q$, the intersection $D^g \cap D'$ is finite.
\item If $\{g_iQg_i^{-1} | 1\leq i\leq n\}$ is a maximal collection of essentially distinct conjugates of $Q$ with infinite intersection. Then $\bigcap g_iQg_i^{-1}$ is a subgroup in $\mc D$ up to conjugacy in $Q$
\end{enumerate}
\end{lemma}
\begin{proof}
By~\cite{GMRS}, there is $n>0$ such that the intersection of any collection of more than $n$ essentially distinct conjugates of $Q$ is finite. By~\cite[Corollary 3.5]{AGM08},  there are only finitely many $Q$-conjugacy classes of subgroups of the form $Q^{g_1}\cap Q^{g_2}\cap \dots \cap Q^{g_j}$ where $j\leq n$ and 
$\{ g_i Q g_i^{-1} : 1\leq i\leq j \}$ is a maximal collection of essentially distinct conjugates with infinite intersection (Here $Q^g$ denotes $gQg^{-1}$). 
 Let $\mc D$ be a finite collection of infinite subgroups obtained by choosing one subgroup of this form for each $Q$-conjugacy class and taking its commensurator in $Q$.  Since  infinite quasiconvex subgroups of  hyperbolic groups have  finite index in their commensurators~\cite{KS96}, $\mc D$ consists of quasiconvex subgroups of $Q$. The other three listed properties follow immediately from the construction.
\end{proof}

\begin{proof}[Proof of Corollary~\ref{cor:3}]
Let $\mc{D}$ be the collection given by Lemma~\ref{lem:qc-core}. The collection $\mc H$ is obtained from $\mc D$ in two steps.
First change $\mc{D}$ to $\mc{D}_0$ by replacing each element of $\mc{D}$ by its commensurator in $G$. Then eliminate redundant entries of $\mc{D}_0$ to obtain $\mc H \subset \mc{D}_0$ which contains no two subgroups conjugate in $G$. Since all elements of $\mc H$ are commensurators of infinite quasiconvex subgroups, they are quasiconvex subgroups~\cite{KS96}. The properties of $\mc D$ given by Lemma~\ref{lem:qc-core} imply
\begin{enumerate}
\item for each $H \in \mc{H}$ and $g \in G\setminus H$, the intersection $H \cap H^g$ is finite,
\item for each distinct pair $H, H' \in \mc{H}$ and $g \in G$, the intersection $H^g \cap H'$ is finite, 
\item for any $H \in \mc{H}$, the intersection $H\cap Q$ is a finite index subgroup of $H$.
\item each infinite intersection of a maximal collection of essentially distinct conjugates of $Q$ is a subgroup of some $H \in \mc H$ up to conjugacy in $H$.
\end{enumerate}
By a result of Bowditch~\cite[Theorem 7.11]{BO99}, the group $G$ is hyperbolic relative to $\mc{H}$. Since $Q$ is a quasiconvex subgroup of $G$, by Theorem~\ref{thm:main} $H$ is a quasiconvex subgroup of $(G, \mc{H})$.
\end{proof}

\begin{remark}
The definition of relative quasiconvexity in the paper~\cite{AGM08} differs from our  Definition~\ref{def:qc}. The equivalence of these two definitions under the assumption that the ambient group is finitely generated (as is the case in Corollary~\ref{cor:3}) is a result by Manning and the author~\cite[Theorem A.10]{MM09}. 
\end{remark}

\subsection{Proof of Corollary~\ref{cor:1}}

\begin{theorem}\cite[Theorem 1.1, Lemma 5.4]{MP09} \label{thm:combination}
Let $(G, \mc H)$ be relatively hyperbolic and let $X$ be a finite generating set of $G$.
For each quasiconvex subgroup $Q$ of $(G, \mc H)$, and each $P \in \mc H$, there is constant $C > 0$ with the following property. 
If $R$ is a subgroup of $H$ such that
\begin{enumerate}
\item $Q \cap H \leq R$, and
\item $|g|_X \geq C$ for any $g \in R \setminus Q$,
\end{enumerate}
then the natural homomorphism $Q \ast_{Q\cap R} R \longrightarrow G$ is injective with image a quasiconvex subgroup of $(G, \mc H)$.  Moreover, every maximal parabolic subgroup $\langle Q\cup R \rangle$ is conjugate in $\langle Q\cup R \rangle$ to either  a maximal parabolic subgroup of $Q$, or to $R$ .
\end{theorem}

\begin{lemma} \label{lem:max-par-hyp}
Let $(G, \mc H)$ be relatively hyperbolic, let $Q$ be a $\sigma$-quasiconvex subgroup of $(G, \mc H)$, and let $P$ be a hyperbolically embedded subgroup of $(G, \mc H)$.  Then maximal parabolic subgroups of $Q$ in $(G, \mc H\cup \{P\})$ are quasiconvex subgroups in $(G, \mc H)$.
\end{lemma}
\begin{proof}
By Theorem~\ref{thm:qc-rh}, maximal parabolic subgroups of $Q$ in $(G, \mc H \cup \{P\})$ are of the form $Q \cap K^g$ where $K \in \mathcal{H}\cup \{P \}$ and $g \in G$. Since $P$ is hyperbolically embedded into $(G, \mc H)$, Corollary~\ref{thm:hyp-qcx} implies that $P$ is quasiconvex in  $(G, \mc H)$. Corollary~\ref{cor:par-qc} implies that each $H \in \mc H$ is quasiconvex in $(G, \mc H)$. Since relative quasiconvexity is preserved by finite intersections and conjugations (Corollary~\ref{prop:qc-conj} and Theorem~\ref{prop:intersection}), the conclusion of the lemma follows.
\end{proof}

\begin{proof}[Proof of Corollary~\ref{cor:1}]
Since $Q$ is $\sigma$-quasiconvex in $(G, \mc H, X, \dist )$, and $P$ is hyperbolically embedded in $(G, \mc H)$, Theorem~\ref{thm:main} implies that $Q$ is quasiconvex in $(G, \mc H \cup \{P\})$.  

By Theorem~\ref{thm:combination} applied to $(G, \mc H\cup \{P\})$, there is a constant $C=C(\mc H, P, X)$ such that the homomorphism $ Q \ast_{Q\cap R} R \rightarrow  G $ is injective whenever $Q\cap P < R < P$ and $|g|_X \geq C$ for all $g\in R\setminus Q$. This is the first conclusion of the Corollary.

Suppose that $R$ is quasiconvex subgroup of $(G, \mc H)$,  $Q\cap P < R < P$, and $|g|_X\geq C$ for all $g\in R\setminus Q$.
We claim that $\langle Q\cup R \rangle$ is also quasiconvex in $(G, \mc H)$. By Theorem~\ref{thm:combination}, $\langle Q\cup R \rangle$ is quasiconvex in $(G, \mathcal{H} \cup \{P \})$. Therefore, by Theorem~\ref{thm:main}, it is enough to show that $\langle Q\cup R \rangle \cap P^g$ is quasiconvex in $(G, \mc H)$ for every $g \in G$. 

Observe that $\langle Q\cup R \rangle \cap P^g$ is a maximal parabolic subgroup of $\langle Q\cup R \rangle$ in $(G, \mc H\cup \{P\})$.
By Theorem~\ref{thm:combination}, $\langle Q\cup R \rangle \cap P^g$ is either conjugate to a maximal parabolic subgroup of $Q$ in $(G, \mc H\cup \{P\})$,  or is conjugate to $R$. By Corollary~\ref{prop:qc-conj} conjugates of $R$ are quasiconvex in $(G, \mc H)$, and by Lemma~\ref{lem:max-par-hyp} maximal parabolic subgroups of $Q$ in $(G, \mc H\cup \{P\})$ are quasiconvex in $(G, \mc H)$. Therefore $\langle Q\cup R \rangle \cap P^g$ is quasiconvex in $(G, \mc H)$.
\end{proof}

\subsection{Proof of Corollary~\ref{cor:2} }

\begin{proposition}\cite[Theorem 1.1]{AM07}\label{lem:am}
Let $G$ be a hyperbolic group with respect to a collection of subgroups $\mc H$. Let $F$ be a finite subset of $G$. Then there exists a collection of subgroups $\mc H'$ with the following properties:
\begin{enumerate}
\item $G$ is hyperbolic relative to $\mc H'$,
\item $\mc H \subset \mc H'$,
\item every $H \in \mc H' \setminus \mc H$ is finite or elementary, and
\item every element of $F$ is parabolic relative to $\mc H'$.
\end{enumerate}
\end{proposition}

\begin{proposition}\label{prop:am}\cite[Lemma 3.8]{AMO07}~\cite[Lemma 8]{AM07}
Let $G$ be a non-elementary and properly relatively hyperbolic group with respect to a collection of subgroups $\mc H$.
Suppose that $G$ has no non-trivial finite normal subgroups. Then there is a hyperbolic element $h$ such that the subgroup $\langle h \rangle$ is hyperbolically embedded into $(G, \mc H)$.
\end{proposition}

\begin{proof}[Proof of Corollary~\ref{cor:2}]
Let $\mc H'$ be the peripheral structure given by Proposition~\ref{lem:am} for the set $F$. By Proposition~\ref{prop:am} there is a hyperbolic element $h$ of infinite order such that the subgroup $\langle h \rangle$ is hyperbolically embedded into $(G, \mc H')$.

Let $f$ be an element of $F$. Since powers of hyperbolic elements are hyperbolic, the cyclic subgroups $\langle f \rangle$ and $\langle h \rangle$ intersect trivially. By Corollary~\ref{cor:1}, there is an integer $n(f)>0$  such that for any $m\geq n(f)$, the subgroup $\langle f , h^m \rangle$ is isomorphic to $\langle f \rangle \ast \langle h^m \rangle$ and quasiconvex in $(G, \mc H')$.  By Corollary~\ref{cor:elementary} the subgroup $\langle f, h^m \rangle$ is quasiconvex in $(G, \mc H)$.
To conclude, let $g = h^m$ where $m$ is the minimum common multiple of $\{n(f) : f \in F\}$.
\end{proof}

\bibliographystyle{plain}
\bibliography{xbib}

\end{document}